\documentclass[12pt]{article}
\usepackage{amssymb}
\def\OPLUS#1{\rb{-5pt}{\mbox{$^{\dis\oplus}_{\,#1}$}}\vs{-4pt}}

\def \Z{\hbox{$Z\hskip -5.2pt Z$}}
\def\OVER#1#2{{\dis{#1\over#2}}}

\def\sZ{\hbox{$\sc Z\hskip -4.2pt Z$}}
\def \Q{\hbox{$Q\hskip -5pt \vrule height 6pt depth 0pt\hskip 6pt$}}
\def \R{\hbox{$I\hskip -4pt R$}}

\def \C{\hbox{$C\hskip -5pt \vrule height 6pt depth 0pt \hskip 6pt$}}

\def\qed{\hfill \hfill \ifhmode\unskip\nobreak\fi\ifmmode\ifinner
         \else\hskip5pt\fi\fi
 \hbox{\hskip5pt\vrule width4pt height6pt depth1.5pt\hskip 1 pt}}
\def\a{\alpha}
\def\b{\beta}
\def\d{\delta}

\def\g{\gamma}
\def\G{\Gamma}
\def\l{\lambda}
\def\L{\Lambda}
\def\e{\epsilon}

\def\DAYU{>_\a\,}
\def\XIYU{<_{{\ssc\,}\a}\,}

\def\Vir{\hbox{\bf\sl Vir}}
\def\HVir{\hbox{{\bf\sl Vir}[$M$]}}

\def\SA{\hbox{$S\!A$}}
\def\SB{\hbox{$S\!B$}}
\def\SVio{\hbox{{\bf\sl SVir}${\ssc\,}_0$}}
\def\SVii{\hbox{{\bf\sl SVir}${\ssc\,}_1$}}
\def\SVir{\hbox{{\bf\sl SVir}[$M,s$]}}

\def\sc{\scriptstyle}
\def\ssc{\scriptscriptstyle}
\def\dis{\displaystyle}
\def\cl{\centerline}

\def\nl{\newline}
\def\ol{\overline}

\def\wh{\widehat}
\def\rar{\rightarrow}
\def\SUM#1#2{\rb{5pt}{\mbox{${}^{#2}_{{}^{\sc\sum}_{\sc#1}}$}}\vs{-4pt}}

\def\Rla{\Leftrightarrow}

\def\bs{\backslash}

\def\vs{\vspace*}
\def\rb{\raisebox}
\def\VS{\vs{4pt}}

\def\AA{{\cal A}}
\def\DD{{\cal D}}

\def\ni{\noindent}
\def\ptl{\partial}
\def\hi{\hangindent}
\def\ha{\hangafter}

\def\Z{\mathbb{Z}}
\def\sZ{\mathbb{Z}}
\def\Q{\mathbb{Q}}
\def\R{\mathbb{R}}
\def\C{\mathbb{C}}

\def\cc{c}
\def\ool#1{[#1)}
\def\oll#1{(#1]}
\def\DIAG{{
$$
\put(15,-30){\put(-150,0)
 {\put(-65,5){{\small (i)}} \put(-2,-2){{\small$\ssc
o$}}\vector(-2,1){30} \put(-25,20)
{{\small$\sc(k,l)$}}\put(-5,0){{\small$\sc b'_1$}}}
\put(-150,0){\vector(-1,2){15}
\put(0,25){{\small$\sc(k',l')$}}\put(8,9){{\small$\sc b'_2$}}}
}
\put(120,-35){ \put(-185,10){{\small (ii)}}
\put(-34,64){{\small$\sc\L':{\ssc\,}(i,j)+r(k,l)+(t+1)(k',l')$}}
\put(-15,39){{\small$\sc (t+1)b'_2$}} \put(-6,-5){{\small$\sc r$}}
\put(14,-15){{\small$\sc q$}} \put(24,-20){{\small$\sc p$}}
\put(-66,25){{\small$\sc s$}}
\put(-43,7){{\small$\sc s-1$}}
\put(-55,60){{\small$\sc b''_1$}}
\put(-32,23){{\small$\sc tb'_2$}}
\put(-62,38){{\small$\sc tb'_2$}}
\put(-36,40){{\small$\sc b''_2$}}
\put(-125,53){{\small$\sc(i,j)+x(k,l)$}}
\put(20,-10){\vector(-2,1){120}}
\put(30,-15){\line(2,-1){20}}
\put(0,0){\vector(-1,2){30}}
\put(-30,15){\vector(-1,2){11}}
\put(-60,30){\vector(-1,2){10}}
\put(-30,60){\bf\vector(-4,-1){38}}
\put(-30,60){\LARGE\line(-4,-1){38}}
\put(-30,60){\LARGE\line(-4,-1){38}}
\put(-30,60){\vector(-4,-1){100}}
\put(-160,30){{\small$\sc\L'+xb''_1$}}
\put(-30,60){\vector(-3,-4){45}}
\put(-142,-5){{\small$\sc\L'+x(k_1b''_1+l_1b''_2)$}}
\put(-30,60){\bf\vector(-1,-2){11}}
\put(-30,60){\LARGE\line(-1,-2){11}}
\put(-30,60){\vector(-1,-2){33}}
\put(-75,-15){{\small$\sc\L'+xb''_2$}}
}
\put(-30,-80){DIAGRAM 1}
$$
}}

\begin{document}
%
%
%
%
%
\cl{\large\bf          Classification of Harish-Chandra Modules
over} \cl{{\large\bf the Higher Rank
         Virasoro Algebras}\footnote{This
         work is supported by a NSF grant 10171064
of China and two grants ``Excellent Young Teacher Program'' and
``Trans-Century Training Programme Foundation for the Talents''
from Ministry of Education of China.\vs{4pt}\nl\hspace*{4.5ex}{\it
Mathematics Subject Classification (1991):} 17B10, 17B65, 17B66,
17B68}}
%
%
%
\vs{2pt}\par
%
%
\ni\cl{\bf Yucai Su}
\par
%
%
\cl{\small\it Department of Mathematics,
                Shanghai Jiaotong University, Shanghai 200030,
                China
                }
                \cl{\small\it
                Email: ycsu@sjtu.edu.cn
}
%
\par\ \vs{-10pt}\par \ni{\small{\bf
                            Abstract:
} We classify the Harish-Chandra modules over the higher rank
Virasoro and super-Virasoro algebras: It is proved that a
Harish-Chandra module, i.e., an irreducible weight module with
finite weight multiplicities, over a higher rank Virasoro or
super-Virasoro algebra is either a module of the intermediate
series, or a finitely-dense module. As an application, it is also
proved that an indecomposable weight module with finite weight
multiplicities over a generalized Witt algebra is either a
uniformly bounded module (i.e., a module with weight
multiplicities uniformly bounded) with all nonzero weights having
the same multiplicity, or a finitely-dense module, as long as the
generalized Witt algebra satisfies one minor condition.
}
\par\
\vs{-10pt}\par
\ni{\bf
             1. \ Introduction
} \vs{-1pt}\par\ni The notions of the higher rank Virasoro
algebras [13] and the higher rank super-Virasoro algebras [17]
appear as natural generalizations of the well-known Virasoro
algebra [4] and super-Virasoro algebras (the Neveu-Schwarz
superalgebra [12], the Ramond superalgebra [14]). As the universal
central extension of the complex Lie algebra of the linear
differential operators over the circle (the Witt algebra), the
Virasoro algebra \Vir, closely related to Kac-Moody algebras
[5--7], is of much interest to both mathematicians and physicists,
partly due to its relevance to string theory [15] and
2-dimensional conformal field theory [3]. The Virasoro algebra
\Vir\ can be defined as a Lie \vs{-4pt}algebra with basis
$\{L_i,\cc\,|\,i\in\Z\}$ such that
$[L_i,L_j]=(j-i)L_{i+j}+{\OVER{i^3-i}{12}}\d_{i,-j}{\ssc\,}\cc,[L_i,\cc]=0$
\vs{-2pt}for all $i,j\in\Z.$ Let $n$ be a positive integer. Let
$M$ be an $n$-dimensional $\Z$-submodule of $\C$, and let $s\in\C$
such that $2s\in M$. A {\it rank $n$ Virasoro algebra} (or a {\it
higher rank Virasoro algebra} if $n\ge2$) [13] is a complex Lie
algebra \HVir\ with basis $\{L_\mu,\cc\,|\,\mu\in M\}$, such
\vs{-7pt}that
$$
[L_\mu,L_\nu]=(\nu-\mu)L_{\mu+\nu}-{\OVER{\mu^3-\mu}{12}}\d_{\mu,-\nu}\cc,\
\ [L_\mu,\cc]=0,\ \ \forall\, \mu,\nu\in M. \vs{-7pt}\eqno(1.1)$$
A {\it rank $n$ super-Virasoro algebra} (or a {\it higher rank
super-Virasoro algebra} if $n\ge2$) [17] is the Lie superalgebra
$\SVir=\SVio\oplus\SVii$, where $\SVio=\HVir$ has a basis
$\{L_\mu,\cc\,|\,\mu\in M\}$ and $\SVii$ has a basis
$\{G_\eta\,|\,\eta\in s+M\}$, with the commutation relations (1.1)
\vs{-9pt}and
$$
\matrix{ [L_\mu,G_\eta]=(\eta-\OVER{\mu}2)G_{\mu+\eta}, \hfill\cr
[G_\eta,G_\l]=2L_{\eta+\l}-\d_{\eta,-\l}\OVER13(\eta^2-\OVER14)\cc,
\ \ [G_\eta,\cc]=0, \hfill\cr} \ \ \forall\,\mu\in M,\,\eta,\l\in
s+M. \vs{-7pt}$$
\hspace*{4ex}
The irreducible representations of the Virasoro and super-Virasoro
algebras, which play very important roles in the theory of vertex
operator (super)algebras and mathematical physics, are well
developed (see for example [2,6--9,16,18]). However for the higher
rank case, not much has been known except modules of the
intermediate series [17,19,20], or Verma-like modules over the
centerless higher rank Virasoro algebras [10]. \vskip -1pt\par%
 A {\it
weight module} is a module $V$ with weight space
\vs{-9pt}decomposition:
$$
V=\,\vs{-7pt}\OPLUS{\!\!\!\!\l\in\C^2}V_\l,\ V_\l=\{v\in
V\,|\,L_0v=\l_0v,\cc v=hv\},\mbox{ where }\l=(\l_0,h)\in\C^2.
$$ A weight module is {\it quasi-finite} if all {\it weight
spaces} $V_\l$ are finite dimensional. A {\it module of the
intermediate series} is an indecomposable quasi-finite weight
module with all the dimensions of weight spaces (and in super case
all the dimensions of weight spaces of ``even'' or ``odd'' part)
are $\le1$. It is proved [17] that a module of the intermediate
series over a higher rank Virasoro algebra $\HVir$ is $A_{a,b},
A_{a'}, B_{a'}$ or one of their quotients for suitable
$a,b,a'\in\C$, where $A_{a,b},A_{a'},B_{a'}$ all have a basis
$\{x_\mu\,|\,\mu\in M\}$ such that $\cc$ acts trivially
\vs{-4pt}and
$$
\matrix{ \hfill A_{a,b}:&L_\mu x_\nu=(a+\nu+\mu
b)x_{\mu+\nu},\VS\hfill&&\cr A_{a'}:&L_\mu
x_\nu=(\nu+\mu)x_{\mu+\nu},\hfill&
\!\!\!\!\!\!\!\!\!\!\!\!\!\!\!\!\!\nu\ne0,\hfill& \hfill L_\mu
x_0&\!\!\!\!= \mu(\mu+a')x_\mu,\VS\hfill\cr B_{a'}:&L_\mu
x_\nu=\nu x_{\mu+\nu},\hfill&
\!\!\!\!\!\!\!\!\!\!\!\!\!\!\!\!\!\nu\ne-\mu,\hfill& \hfill\ \ \
L_\mu x_{-\mu}&\!\!\!\!= -\mu(\mu+a')x_0,\hfill\cr }
\eqno\matrix{(1.2{\rm a})\VS\cr(1.2{\rm b})\VS\cr(1.2{\rm
c})\cr}\vs{-4pt}$$ for $\mu,\nu\in M$; and that a module of the
intermediate series over the higher rank super-Virasoro algebras
$\SVir$ is one of the three series of the modules
$\SA_{a,b},\SA_{a'},\SB_{a'}$ or their quotient modules for
suitable $a,b,a'\in\C$, where $\SA_{a,b},\SA_{a'}$ have basis
$\{x_\mu\,|\,\mu\in M\}\cup\{y_\eta\,|\,\eta\in s+M\}$ and
$\SB_{a'}$ has basis $\{x_\eta\,|\,\eta\in
s+M\}\cup\{y_\mu\,|\,\mu\in M\}$ such that $\cc$ acts trivially
\vs{-3pt}and
$$
\matrix{
\matrix{
\SA_{a,b}\!:\!\!\!\!\hfill&
L_\mu x_\nu=(a+\nu+\mu b)x_{\mu+\nu},\hfill&
L_\mu y_\eta=(a+\eta+\mu(b-\OVER12))y_{\mu+\eta},
\hfill\cr&
G_\l x_\nu=y_{\l+\nu},\hfill&
G_\l y_\eta=(a+\eta+2\l(b-\OVER12))x_{\l+\eta},
\hfill\cr
}\vs{5pt}\hfill\cr\matrix{
\SA_{a'}\!:\!\!\!\!\hfill&
L_\mu x_\nu=(\nu+\mu)x_{\mu+\nu},\,\nu\ne0,\hfill&
L_\mu x_0=\mu(\mu+a')x_\mu,\hfill&
L_\mu y_\eta=(\eta+\OVER{\mu}2)y_{\mu+\eta},
\vs{2pt}\hfill\cr&
G_\l x_\nu=y_{\l+\nu},\,\nu\ne0,\hfill&
G_\l x_0=(2\l+a')y_\l,\hfill&
G_\l y_\eta=(\eta+\l)x_{\l+\eta},
\hfill\cr
}\vs{5pt}\hfill\cr\matrix{
\SB_{a'}\!:\!\!\hfill&
L_\mu x_\eta=(\eta+\OVER{\mu}2)x_{\mu+\eta},\hfill&
L_\mu y_\nu=\nu y_{\mu+\nu},\,\nu\ne-\mu,\hfill&
L_\mu y_{-\mu}=-\mu(\mu+a')y_0,
\vs{2pt}\hfill\cr&
G_\l x_\eta=y_{\l+\eta},\,\eta\ne-\l,\hfill&
G_\l x_{-\l}=(2\l+a')y_0,\hfill&
G_\l y_\nu=\nu x_{\l+\nu},
\hfill\cr}
\hfill\cr}
\vs{-4pt}$$
for $\mu,\nu\in M,\,\l,\eta\in s+M$.
\par
A {\it uniformly bounded module} is a quasi-finite weight module
with all weight multiplicities being uniformly bounded. It is
proved [19,20] that a uniformly bounded irreducible module is
simply a module of the intermediate series. An irreducible
quasi-finite weight module is called a {\it Harish-Chandra
module}. It is a deep fact (first conjectured in [6], then proved
in [9] and partially proved in [2,8,16]) that a Harish-Chandra
module over the Virasoro algebra is either a module of the
intermediate series or else a highest or lowest weight module. The
same is true for the super-Virasoro algebras [18]. Unlike the
Virasoro or super-Virasoro case, a nontrivial highest or lowest
weight module, or more generally, a nontrivial Verma-like module
[10] for a higher rank Virasoro or super-Virasoro algebra
is not a quasi-finite weight module. 
\vs{-1pt}\par The classification of Harish-Chandra modules
is definitely an important problem in the representation theory of
Lie algebras. Since no Harish-Chandra modules other than modules
of the intermediate series have been found for the higher rank
Virasoro or super-Virasoro algebras, a natural question is:
\vs{-1pt}\par {\it Does there exist any Harish-Chandra module
other than modules of the intermediate series over the higher rank
Virasoro or super-Virasoro algebras?} \vs{-1pt}\par The aim of
this paper is to answer this question, mainly we have the
following theorem (for the $\Q$-Virasoro algebra, this theorem was
obtained in [11]).
\par\ni
{\bf Theorem 1.1.} {\it A Harish-Chandra module over a higher rank
Virasoro or super-Virasoro algebra is either a module of the
intermediate series, or a finitely-dense module (cf.~(2.10)).}
\par
Let $m$ be a positive integer. Let $\G$ be a nondegenerate
additive subgroup of $\C^m$, that is, $\G$ contains a $\C$-basis
of $\C^m$. Let $\pi_i:\C^m\rar\C$ be the natural projection from
an element of $\C^m$ to its $i$th coordinate, $i=1,2,...,m$. Let
$\AA=\C[\G]$ be the group algebra of $\G$ with basis
$\{x^\a\,|\,\a\in\G\}$. Define the derivations $\ptl_i$ of $\AA$
such that $\ptl_i(x^\a)=\pi_i(\a)x^\a$. Set $\DD={\rm
span}\{\ptl_i\,|\,i=1,2,...,m\}$. Then the tensor space
$\AA\otimes\DD$ can be defined as a Lie algebra under the bracket
$[x^\a\ptl_i,x^\b\ptl_j]=x^{\a+\b}(\pi_i(\b)\ptl_j-\pi_j(\a)\ptl_i)$,
called a {\it generalized Witt algebra}, denoted by $W(m,\G)$ (see
for example [23]). Note that a generalized Witt algebra defined in
[10] is simply a Lie algebra $W(1,\G)$ with $\G\subset\C$ being a
free $\Z$-module of finite rank, i.e., it is a centerless higher
rank Virasoro algebra using our notion (1.1)). Observe that a
generalized Witt algebra $W(m,\G)$ has a (unique) nontrivial
central extension if and only if $m=1$, in this case the universal
central extension of \vs{-3pt}$W(1,\G)$ is referred to as a {\it
generalized Virasoro algebra} in [22], denoted by $\wh W(1,\G)$.
As an application of Theorem 1.1, we obtain
\par\ni
{\bf Theorem 1.2.} {\it Suppose there is a group injection
$\Z\times\Z\rar\G$ (in particular, if $m\ge2$ there always exists
a group injection $\Z\times\Z\rar\G$ by \vs{-2pt}nondegenerateness
of $\G$). Then an indecomposable quasi-finite weight module over
$W(m,\G)$ or $\wh W(1,\G)$ is either a uniformly bounded module
with all nonzero weights having the same multiplicity, or a
finitely-dense module.}
\par We believe that Theorem 1.2 will certainly be important to the
classification problem of the Harish-Chandra modules over the
generalized Witt algebras (especially the classical Witt algebras
$W_n={\rm Der\,}\C[t_1^\pm,...,t_m^\pm]$, the derivation algebras
of the Laurent polynomial algebras of $m$ variables, which is
closely related to toroidal Lie algebras, cf.~[1]). This is also
one of our motivations to present results here. \vs{-1pt}\par
The paper is organized as follows. In Section 2, after collecting
some necessary information on representations of the higher rank
Virasoro algebras, we give the proof of Theorem 1.1 for the case
$n=2$. Then we complete the proofs of Theorems 1.1 and 1.2 in
Section 3.
\par\
\vs{-10pt}
\par\ni{\bf 2. \ Preliminaries and the Case $n=2$}
\vs{-1pt}\par\ni First observe that by regarding a module over a
higher rank super-Virasoro algebra as a module over a higher rank
Virasoro algebra, using Theorem 1.2, one obtains that a
not-finitely-dense Harish-Chandra module over a higher rank
super-Virasoro algebra is uniformly bounded, and thus using
results in [20], it is a module of the intermediate series. Thus
in the following, we shall only consider the ``unsuper'' case.
\par
So, let $\HVir$ be a higher rank Virasoro algebra, where $M$ is an
$n$-dimensional $\Z$-submodule of $\C$ with $n\ge2$. Following
[19], a \HVir-module $V$ is called a {\it GHW module} (here GHW
stands for ``generalized highest weight'') if $V$ is generated by
a nonzero weight vector $v_{\ssc\L}$ and there exists a $\Z$-basis
$B=\{b_1,...,b_n\}$ of $M$ such \vs{-4pt}that
$$
L_\mu v_{\ssc\L}=0\mbox{ for all }
0\ne\mu=\SUM{i=1}{\ n}m_i b_i\in M
\mbox{ with all coefficients }m_i\in\Z_+.
\eqno(2.1)$$
The vector $v_{\ssc\L}$ is called a {\it GHW vector with respect to
$\Z$-basis $B$} or simply a {\it GHW vector}.
\par
Highest or lowest weight modules, no matter how the ordering on
$M$ is defined [11,13], are GHW modules. The following three
lemmas are taken from [19], [21], [9] respectively.
\par\ni{\bf Lemma 2.1.} {\it (1) A uniformly bounded Harish-Chandra
$\Vir[M]$-module is a module of the intermediate series. (2) A
Harish-Chandra $\Vir[M]$-module is either a module of the
intermediate series or a GHW module.}
\par\ni
{\bf Lemma 2.2.}
{\it Let $W$ be a quasi-finite weight $\Vir$-module such
that $\l=(\l_0,h)$ is a weight.
For any $i\in\Z$, there exist only a finite number of primitive vectors with
weights $\l+j=(\l_0+j,h)$ such that $j\ge i$ (a primitive vector is a nonzero
weight vector $v$ such that $L_iv=0$ for $i\ge1$).}
\par\ni
{\bf Lemma 2.3.}
{\it Let $W$ be a quasi-finite weight $\Vir$-module.
Suppose $0\ne v\in W$ has weight $\l$ such that $L_iv=0$ for $i>>0$. Then
there exists a primitive vector with weight $\l+j$ for some $j\ge0$.}
\par
By Lemma 2.1, the proof of Theorem 1.1 is equivalent to proving
that a nontrivial not-finite-dense quasi-finite irreducible GHW
module does not exist. In this section, we shall prove Theorem 1.1
for the special case $n=2$.
\par
Thus suppose $V$ is a nontrivial not-finite-dense quasi-finite
irreducible GHW $\HVir$-module generated by a GHW vector
$v_{\ssc\L}$ with weight $\L=(\L_0,h)\in\C^2$, where $M=\Z b_1+\Z
b_2$ such that $b_1,b_2\in\C$ are $\Z$-linear independent and
$$
L_{i,j}v=0\mbox{ for all }(i,j)\in\Z^2\mbox{ with }0\ne(i,j)\ge0,
\mbox{ where }L_{i,j}=L_{ib_1+jb_2}. \eqno(2.2)$$ Here in general,
we define $(i,j)\ge(k,l)$ if $i\ge k,j\ge l$ and define
$(i,j)>(k,l)$ if $i>k,j>l$.
For $m\in M,\,S\subset M$, we denote by $\Vir[m]$ or $\Vir[S]$ the
subalgebra of $\HVir$ generated by $\{L_{\pm m},L_{\pm2m}\}$ or
$\{L_{\pm m},L_{\pm2m}\,|\,m\in S\}$ respectively. In particular,
if $m\in M\bs\{0\}$, then $\Vir[m]$ is a rank one Virasoro
subalgebra isomorphic to $\Vir$. If $\emptyset\ne S\subset
M\bs\{0\}$, then $\Vir[S]$ is a rank one or rank two Virasoro
subalgebra of $\HVir$. For any subalgebra $L$ of $\HVir$, we
denote by $U(L)$ the universal enveloping algebra of $L$.
\par
Since for an irreducible weight module, the central element $\cc$
must act as a scalar $h$ for some fixed $h\in\C$, we shall always
omit $h$ and simply denote a weight $\l=(\l_0,h)$ by
$\l=\l_0\in\C$. Let $P$ be the set of weights of $V$. Then
$P\subset\L+\Z b_1+\Z b_2$. Thus when $b_1,b_2$ are fixed, we can
define an injection
$$
\phi_{b_1,b_2}:P\rar\Z^2\mbox{ \ such that \ }
\phi_{b_1,b_2}(\L+ib_1+jb_2)=(i,j).
\eqno(2.3)$$
We extend $\phi_{b_1,b_2}$ to $\phi_{b_1,b_2}:\L+\Z b_1+\Z b_2\rar\Z^2$.
To better understand the following discussion,
one can regard $\Z^2$ as the set of integral points in the $Oxy$-plane.
\par
For $m\in\Z$, we denote $\oll{-\infty,m}=\{i\in\Z\,|\,i\le m\}$
and $\ool{m,\infty}=\{i\in\Z\,|\,i\ge m\}$.
\par\ni
{\bf Lemma 2.4.}
{\it For any $v\in V$, there exists $p>>0$ such that
$L_{i,j}v=0$ for all $(i,j)\ge(p,p)$.
}
\par\ni
{\it Proof.} Since $v=uv_{\ssc\L}$ for some $u\in U(\HVir)$, then
$u$ is a linear combination of elements of the form
$u'=L_{i_1,j_1}\cdots L_{i_n,j_n}$. Thus without loss of
generality, we can suppose $u=u'$. Take $p={\rm
max}\{-\sum_{i_q<0}i_q,-\sum_{j_q<0}j_q\}+1$. Using relation
$[u,u_1u_2]=[u,u_1]u_2+u_1[u,u_2]$ for $u\in\HVir$ and $u_1,u_2\in
U(\HVir)$, and using relations (1.1) and by (2.2), we obtain
$L_{i,j}v=[L_{i,j},u]v_{\ssc\L}=0$ for all $(i,j)\ge(p,p)$ by
induction on $n$. \qed\par For any $p\in\Z$, it is easy to see
that
$$
B'=\{b^*_1=(p+1)b_1+(p+2)b_2,\ b^*_2=pb_1+(p+1)b_2\}, \eqno(2.4)$$
is also a $\Z$-basis of $M$. Thus by Lemma 2.4 and definition
(2.1), by choosing $p>>0$, one sees that every nonzero weight
vector is a GHW vector with respect to the $\Z$-basis $B'$.\par\ni
{\bf Lemma 2.5.} {\it Let $\Vir[S]$ be any rank two Virasoro
subalgebra of $\Vir[M]$ and let $b^*_1,b^*_2$ be any $\Z$-basis of
$M$. Then $\Vir[M]$ is generated by
$\Vir[S]\cup\{L_{b^*_1},L_{b^*_2}\}$.}
\par\ni
{\it Proof.} Since $\Vir[S]$ is a rank two Virasoro subalgebra of
$\HVir$, there exists $(i,j)>0$ such that
$L_{-ib^*_1-jb^*_2}\in\Vir[S]$. Clearly $L_{-b^*_1},L_{-b^*_2}$
can be generated by $L_{-ib^*_1-jb^*_2}, L_{b^*_1},L_{b^*_2}$.
Since $\HVir$ is generated by $L_{\pm b^*_1},L_{\pm b^*_2}$, it is
also generated by $\Vir[S]\cup\{L_{b^*_1},L_{b^*_2}\}$.
\qed\par\ni {\bf Lemma 2.6.} {\it Let $\Vir[S]$ be a rank two
Virasoro subalgebra of $\Vir[M]$. Then any nonzero
$\Vir[S]$-submodule of $V$ is nontrivial.
}\par\ni {\it Proof.} Suppose $\C v\ne\{0\}$ is a trivial
$\Vir[S]$-submodule. Choose $p>>0$ satisfying Lemma 2.4 and let
$b^*_1,b^*_2$ be as in (2.4). Then
$\Vir[S]\cup\{L_{b^*_1},L_{b^*_2}\}$ acts trivially on $\C v$.
Since $V$ is irreducible, by Lemma 2.5, $V=\C v$ is a trivial
$\HVir$-module, a contradiction with the assumption that $V$ is
nontrivial. \qed\par\ni {\bf Lemma 2.7.} {\it For any $v\in
V\bs\{0\}$ and any $(i,j)\in\Z^2$ with $(i,j)>0$, we have
$L_{-i,-j}v\ne0.$ }
\par\ni
{\it Proof.} Suppose $L_{-i,-j}v=0$ for some $(i,j)>0$. Let $p$ be
as in Lemma 2.4, then $L_{-i,-j}$, $L_{pi+1,pj},L_{pi,pj+1}$ act
trivially on $v$, but $\HVir$ is generated by these three
elements, a contradiction with that $V$ is a nontrivial
irreducible module. \qed\par\ni {\bf Lemma 2.8.} {\it Let
$(i,j)\in\phi_{b_1,b_2}(P)$ and let $(k,l)>0$. Then
$$\{x\in\Z\,|\,(i,j)+x(k,l)\in\phi_{b_1,b_2}(P)\}=\oll{-\infty,m}
\;\;\;\;\mbox{for \,some}\;\;\;\;m\ge0.\eqno(2.5)$$}
\vs{-4pt}\par\ni {\it Proof.} Denote by $I$ the set in the
left-hand side of (2.5). By Lemma 2.7, if $x_1\in I$, then
$\oll{-\infty,x_1}\subset I$. Suppose $I=\Z$. For any $x_2\in\Z$,
fix a weight vector $v_{x_2}$ of weight
$\phi^{-1}_{b_1,b_2}((i,j)+x_2(k,l))$. By Lemma 2.4, when $x>>0$,
$L'_x v_{x_2}=0$, where $L'_m=L_{mk,ml},m\in\Z$ span a rank one
Virasoro subalgebra $\Vir[kb_1+lb_2]$. Thus by Lemma 2.3, there
exists $x_3\ge x_2$ such that
$\phi^{-1}_{b_1,b_2}((i,j)+x_3(k,l))$ is a
$\Vir[kb_1+lb_2]$-primitive weight. This contradicts Lemma
2.2.\qed\par By Lemma 2.8, we can suppose
$\{x\in\Z\,|\,x(1,1)\in\phi_{b_1,b_2}(P)\} =\oll{-\infty,p-2}$ for
some $p\ge2$. Then $(p-1,p-1)\notin\phi_{b_1,b_2}(P)$ and Lemma
2.8 again shows that $(i,j)\notin \phi_{b_1,b_2}(P)$ for all
$(i,j)\ge(p,p)$ which is $>(p-1,p-1)$. Taking $b^*_1,b^*_2$ to be
as in (2.4), by (2.3), we have
$$
\matrix{
\phi_{b^*_1,b^*_2}^{-1}(i,j)
\!\!\!\!&=\L+ib^*_1+jb^*_2
\vs{4pt}\hfill\cr&
=\L+(i(p+1)+jp)b_1+(i(p+2)+j(p+1))b_2
\vs{4pt}\hfill\cr&
=\phi^{-1}_{b_1,b_2}(i(p+1)+jp,i(p+2)+j(p+1))
\mbox{ for }(i,j)\in\Z^2.
\hfill\cr}
$$
Thus by Lemmas 2.7 and 2.8, for $(i,j)\in\Z^2$, we have
$$
\matrix{ (i,j)\notin\phi_{b^*_1,b^*_2}(P)\mbox{ if }0\ne(i,j)\ge0,
\mbox{ and } (i,j)\in\phi_{b^*_1,b^*_2}(P)\mbox{ if }(i,j)\le0,
\mbox{ and} \vs{7pt}\hfill\cr
(k,l)\notin\phi_{b^*_1,b^*_2}(P)\mbox{ for all }(k,l)\ge(i,j)
\mbox{ if }(i,j)\notin\phi_{b^*_1,b^*_2}(P), \mbox{ and }
\vs{7pt}\hfill\cr (k,l)\in\phi_{b^*_1,b^*_2}(P)\mbox{ for all
}(k,l)\le(i,j) \mbox{ if }(i,j)\in\phi_{b^*_1,b^*_2}(P),
\hfill\cr} \eqno\matrix{(2.6)\vs{7pt}\cr\ \vs{7pt}\cr(2.7)\cr}$$
and for any $(i,j)\in\phi_{b^*_1,b^*_2}(P)$ and any $(k,l)\ge0$,
$$
\{x\in\Z\,|\,(i,j)+x(k,l)\in\phi_{b^*_1,b^*_2}(P)\}=\oll{-\infty,m}
\mbox{ for some }m\ge0. \eqno(2.8)$$ In the following, for
convenience, we redenote $b^*_1,b^*_2$ by $b_1,b_2$ respectively
and from now on we simply denote $\phi_{b_1,b_2}$ by $\phi$.
\par
The following rather technical lemma is crucial in obtaining our
classification.
\par\ni
{\bf Lemma 2.9.} {\it For any $(i,j),(k,l)\in\Z^2$, there exists
$p\in\Z$ such that
$$
\{x\in\Z\,|\,(i,j)+x(k,l)\in\phi(P)\}=\oll{-\infty,p}\mbox{ \ or \
}\ool{p,\infty}.$$%
Proof.} %
Let $I=\{x\in\Z\,|\,(i,j)+x(k,l)\in\phi(P)\}$. Suppose
$(k,l)=d(k_1,l_1)$ for some $d,k_1,l_1\in\Z$. Clearly if the
result holds for $(k_1,l_1)$, then it holds for $(k,l)$ too. Thus
by replacing $(k,l)$ by $(k_1,l_1)$ we can suppose $k,l$ are
coprime and $k\le0$. Suppose $I\ne\oll{-\infty,p}$ and
$I\ne\ool{p,\infty}$ for any $p$. Then by (2.8) we must have
$k<0<l$.
\par
{\it Case 1:} $I=\emptyset.$ Since
$(0,0)\in\phi(P)$, we obtain $(i,j)\notin\Z(k,l)$.
Furthermore, since $k,l$ are coprime, there do not exist $p,q\in\Z\bs\{0\}$
such that $p(k,l)+q(i,j)=0$, i.e., $(i,j), (k,l)$ are $\Z$-linear
independent. So $\Vir[b'_1,b'_2]$ is a rank two Virasoro subalgebra
of $\HVir$, where $b'_1=kb_1+lb_2,b'_2=ib_1+jb_2$.
Then $I=\emptyset$ means that
$$
\L+nb'_1+b'_2\mbox{ is not a weight for any }n\in\Z.
\eqno(2.9)$$
Let $W=\oplus_{\l\in\L+\sZ b'_1}V_\l$, which
is a $\Vir[b'_1]$-submodule of $V$.
\par
{\bf Claim 1.} $W$ is a uniformly bounded $\Vir[b'_1]$-module.
\par
Let $0\ne w\in V_{\L+mb'_1}$ for some $m\in\Z$. If
$L_{-mb'_1-b'_2}w=0$, then $\C w$ is a trivial
$\Vir[b'_1,b'_2]$-submodule since $\Vir[b'_1,b'_2]$ can be
generated by the set
$\{L_{-mb'_1-b'_2},L_{nb'_1+b'_2}\,|\,n\in\Z\}$ which acts
trivially on $w$ by (2.9). This is contrary to Lemma 2.6. Thus we
obtain a linear injection
$L_{-mb'_1-b'_2}|_{V_{\L+mb'_1}}:V_{\L+mb'_1}\rar V_{\L-b'_2}$.
Thus ${\rm dim\ssc\,}V_{\L+mb'_1} \le{\rm dim\ssc\,}V_{\L-b'_2}$
for all $m\in\Z$, i.e., $W$ is uniformly bounded and the claim is
proved.
\par
Since a uniformly bounded module has only a finite number of
composition factors, we can take an irreducible
$\Vir[b'_1]$-submodule $W'$ of $W$. Let $W''$ be the
$\Vir[b'_1,b'_2]$-submodule of $V$ generated by $W'$. Then $W''$
is a quotient module of a finitely-dense generalized Verma module
[10], thus $V$ is a finitely-dense module in the following sense.
So this case does not occur.
\par\ni%
{\bf Definition 2.10}. A module $V$ over a Lie algebra $\cal L$ is
a {\it finitely-dense module} if there exist a rank 2 Virasoro
subalgebra ${\cal L}'$ of $\cal L$ and a ${\cal L}'$-submodule
$V'$ of $V$ such that $V'$ is a quotient of a finitely-dense
generalized Verma module. A module $V$ over a rank 2 Virasoro
algebra $\Vir[M]$ is a {\it finitely-dense generalized Verma
module} if there exist a $\Z$-basis $\{b_1,b_2\}$ of $M$ and an
intermediate series $\Vir[b_1]$-submodule $V'$ of $V$ such that
$$%
V={\rm Ind}_{\textsl{\footnotesize
Vir}[b_1]\oplus{\textsl{\footnotesize
Vir}}^+[M]}^{{\textsl{\footnotesize
Vir}}[M]}V'=U(\Vir[M])\otimes_{U(\textsl{\footnotesize
Vir}[b_1]\oplus{\textsl{\footnotesize Vir}}^+[M])}V'\cong
U(\Vir^-[M])\otimes_{\C}V',\eqno(2.10)
$$%
where $\Vir^{\pm}[M]={\rm span}\{L_{mb_1+nb_2}\,|\,m\in\Z,\,\pm
n\in\Z_+\bs\{0\}\}$.
\par
{\it Case 2:}
There exist $x_1,x_2,x_3\in\Z$, $-x_1<x_2<x_3$ such that
$-x_1,x_3\in I$ but $x_2\notin I$, i.e.,
$$
(i,j)-x_1(k,l)\notin\phi(P),\ (i,j)+x_2(k,l)\in\phi(P),\
(i,j)+x_3(k,l)\notin\phi(P). \eqno(2.11)$$ Replacing $x_2$ by the
largest $x<x_3$ such that $(i,j)+x(k,l)\in\phi(P)$, and then
replacing $x_3$ by $x_2+1$ and $(i,j)$ by $(i,j)+x_2(k,l)$, we can
suppose $-x_1<x_2=0,x_3=1$. Choose a nonzero weight vector $v_\l$
with weight $\l=\phi^{-1}(i,j)$. Then (2.11) means that
$L_{k,l}v_\l=0=L_{-x_1k,-x_1l}v_\l$, thus also $L_{-k,-l}v_\l=0$.
Choose $p,q>>0$ such that $L_{p,q}v_\l=0$ (cf.~Lemma 2.4). Then
$S=\{b'_1=kb_1+lb_2,b'_2=pb_1+qb_2\}$ is a $\Z$-linear independent
subset of $M$ since $kq-lp<0$ (note that $k<0<l$). Note that for
any $mb'_1+nb'_2\in M$, if $n>0$, then $L_{mb'_1+nb'_2}$ can be
generated by $L_{\pm b'_1},L_{b'_2}$. Thus $L_{mb'_1+nb'_2}v_\l=0$
if $n>0$. For any $a\in U(\Vir[S])$ with weight $mb'_1+nb'_2$ such
that $n>0$, we can write $a$ as a linear combination of the form
$$
L_{i_1b'_1+j_1b'_2}\cdots L_{i_sb'_1+j_sb'_2}\mbox{ with }j_s>0.
\vs{-2pt}$$
Thus $av_\l=0$. Let $V'=U(\Vir[S])v_\l$, this shows that
$\l+mb'_1+nb'_2$ is not a weight of $V'$ for any $n>0$.
Now as discussion in Case 1, we obtain that $V'$ is not quasi-finite.
Thus this case does not occur either.
\par
{\it Case 3:} $I=\oll{-\infty,p}\cup\ool{q,\infty}$ for some
$p,q\in\Z$. Since $k,l$ are coprime and $k<0<l$, we can choose
$k',l'$ such that
$$kl'-lk'=-1\;\;\;\;\mbox{and}\;\;\;\;k'<0<l'.\eqno(2.12)$$
Then
$$
\{b'_1=kb_1+lb_2,b'_2=k'b_1+l'b_2\},$$ forms a
$\Z$-basis of $M$. Then $b'_1,b'_2$ are shown as in Diagram 1(i).
\par
\DIAG \vs{4pt}\par Assume that there exists $t>0$ such that
$(i,j)+r(k,l)+(t+1)(k',l')\notin \phi(P)$ for all $r\in I$. This
means that $$\{x\in\Z\,|\,(i',j')+x(k,l)\in\phi(P)\} \subset\Z\bs
I\;\;\;\;\mbox{for}\;\;\;\;(i',j')=(i+(t+1)k',j+(t+1)l').$$
 But $\Z\bs I$
is either empty or a finite set $\{p+1,...,q-1\}$. This falls to
Case 1 or Case 2, which is impossible. Thus, we obtain that for
any $t>0$, there always exists $r\in I$ such that
$$(i,j)+r(k,l)+(t+1)(k',l')\in\phi(P).\eqno(2.13)$$
Fix $x\in\Z$ such that $x>>0$, and let $s=(k'-l')x,t=(l-k)x$.
Since $k<0<l, k'<0<l'$, we have $t>>0,s<<0$. In particular, we can
assume $s<p$, so $s\in I$. Since $(m,n)\notin\phi(P)$ for all
$(m,n)>0$, and by (2.12) we have $s(k,l)+t(k',l')=x(1,1)>>0$, thus
we can assume
$$
(i,j)+s(k,l)+t(k',l'),\ (i,j)+(s-1)(k,l)+t(k',l')\notin\phi(P).
\eqno(2.14)$$ This proves that we can choose $r,s\in I,t\ge0$ such
that (2.13), (2.14) hold.
Say $r<s-1$ (otherwise the proof is similar). 
Now take
$$
b''_1=(s-r)b'_1-b'_2,\ \ b''_2=(s-r-1)b'_1-b'_2.$$
Then $\{b''_1,b''_2\}$ is still a $\Z$-basis of $M$ (cf.~Diagram
1(ii)). By (2.13), we can choose a nonzero weight vector
$v_{\ssc\L'}$ with the weight $\L'=\phi^{-1}
((i,j)+r(k,l)+(t+1)(k',l'))$. Then by (2.14),
$$
L_{b''_1}v_{\ssc\L'}=L_{b''_2}v_{\ssc\L'}=0,\mbox{ and thus }
L_{mb''_1+nb''_2}v_{\ssc\L'}=0\mbox{ for all }(m,n)>0.
$$
Using this, one sees that $v_{\ssc\L'}$ is a GHW vector with
respect to the $\Z$-basis $\{b''_1+b''_2,b''_1+2b''_2\}$ of $M$,
thus as in the proofs in Lemmas 2.4, 2.8, for any $v\in V$, we
have $L_{mb''_1+nb''_2}v=0$ for all $(m,n)>>0,$ and for any
$(k_1,l_1)>0$, we have
$$
\L'+x(k_1b''_1+l_1b''_2)\mbox{ is not a weight for }x>>0.
\eqno(2.15)$$ From Diagram 1(ii), one can easily see that when
$x>>0$, under mapping (2.3), $\L'+x(k_1b''_1+l_1b''_2)$
corresponds to a point located below and to the left of some point
$(i,j)+y(k,l)$ for some $y\in I$, which by (2.7) is in $\phi(P)$,
contrary to (2.15). Thus this case does not occur either. This
completes the proof of Lemma 2.9.
\qed\par\ni%
{\it Proof of Theorem 1.1 for the case $n=2$.} For any
$i\in\Z_+\bs\{0\}$, since $(-i,-j)\in\phi(P)$ for all $j\in\Z_+$
(cf.~Lemma 2.7), by Lemma 2.9, we can \vs{-3pt}let
$$
y_i={\rm max}\{y\in\Z\,|\,(-i,y)\in\phi(P)\},\
x_i={\rm max}\{x\in\Z\,|\,(x,-i)\in\phi(P)\}.
\vs{-3pt}\eqno(2.16)$$
By (2.7), we have $y_{i+1}\ge y_i\ge0,x_{i+1}\ge x_i\ge0$.
For $y\in\Z$, if $y\le y_i$ then
since $(0,0),(-i,y)\in\phi(P)$, by Lemma 2.9, all integral points
lying between
$(0,0)$ and $(-i,y)$ and
lying on the line linking them
must be in $\phi(P)$, i.e.,
$$
\Z^2\cap\{x(-i,y)\,|\,x\in\Q,0\le x\le1\}\subset\phi(P) \mbox{ for
}i\in\Z_+\bs\{0\},\,y\in\Z_+,\,y\le y_i. \vs{-1pt}\eqno(2.17)$$
Let $j,t\in\Z_+\bs\{0\}$. If $y_{tj}\ge t(y_j+1)$, then
$t(-j,y_j+1)=(-tj,t(y_j+1))\in\phi(P)$, and by (2.17) it gives
$(-j,y_j+1)\in\phi(P)$, contrary to definition (2.16). This
\vs{-3pt}means
$$
y_{tj}<t(y_j+1) 
\mbox{ for all }t,j\in\Z_+\bs\{0\}.
\vs{-3pt}\eqno(2.18)$$
Since $(0,1)\notin\phi(P)$ (cf.~(2.6)), but $(-j,y_j)\in\phi(P)$,
by Lemma 2.9, we also \vs{-3pt}have
$$
\Z^2\cap\{(0,1)+x(-j,y_j-1)\,|\,x\in\Q,x\ge1\}\subset\phi(P)
\mbox{ for }j\in\Z_+\bs\{0\}.
\vs{-2pt}\eqno(2.19)$$
In particular, taking $x=t\in\Z_+\bs\{0\}$, it implies $(-tj,t(y_j-1)+1)\in
\phi(P)$, \vs{-2pt}i.e.,
$$
y_{tj}\ge t(y_j-1)+1
\mbox{ for }
t,j\in\Z_+\bs\{0\}.
\vs{-1pt}\eqno(2.20)$$
Using (2.18), (2.20), we obtain
$j(y_i\!-\!1)\!+\!1\!\le\! y_{ij}\!<\!i(y_j\!+\!1)$
for $i,j\!\in\!\Z_+\bs\{0\}$. From this we
\vs{-3pt}deduce
$$
\OVER{y_j}{j}-\OVER{i+j-1}{ij}<\OVER{y_i}{i}<\OVER{y_j}{j}+\OVER{i+j-1}{ij}
\mbox{ \ for all }i,j\in\Z_+\bs\{0\}.
\vs{-3pt}$$
This shows that the following limits
\vs{-5pt}exist:
$$
\a=\lim_{i\rar\infty}\OVER{y_i}{i},\
\b=\lim_{i\rar\infty}\OVER{x_i}{i}, \vs{-5pt}\eqno(2.21)$$ where
the second equation is obtained by symmetry. Since
$(0,1)\notin\phi(P)$ and $(-j,y_j)\in\phi(P)$, we have
$(-xj,x(y_j-1)+1)=(0,1)+x(-j,y_j-1) \notin\phi(P)$ for all
$x\in\oll{-\infty,0}$ (cf.~(2.19)). In particular, taking $x=-t$
and by definition (2.16), it \vs{-3pt}means
$$
x_{t(y_j-1)-1}<tj\mbox{ for all }t,j\in\Z_+\bs\{0\}\mbox{ such that }
t(y_j-1)-1\ge1.
\vs{-7pt}\eqno(2.22)$$
Note that Lemma 2.9 means that for any $i\in\Z_+\bs\{0\}$
there exists $j\in\Z_+\bs\{0\}$ such that $(-j,i)\in\phi(P)$,
i.e., $\lim_{j\rar\infty}y_j=\infty$. Dividing (2.20) by
$tj$ and taking $\lim_{{\ssc\,}t\rar\infty}$, we obtain $\a>0$.
Dividing (2.22) by $t(y_j-1)-1$ and taking $\lim_{j\rar\infty}$, we obtain
$\b\le\a^{-1}$.
\vs{-1pt}\par
For $j\in\Z_+\bs\{0\}$, we \vs{-5pt}denote
$$
I_j=\{x\in\Q\,|\,x\ge0, -(1,1)+x((-j,y_j)+(1,1))\in\phi(P)\}.
\vs{-3pt}$$
Taking $x=j+1$, we have
$$
-(1,1)+(j+1)((-j,y_j)+(1,1))=(-j^2,-1+(j+1)(y_j+1))\notin\phi(P),
$$
because $-1+(j+1)(y_j+1)\ge j(y_j+1)>y_{j^2}$ by (2.18). Thus
$j+1\notin I_j$. So $I_j$ is a finite set, and then Lemma 2.9
implies
$$
\Z^2\cap\{-(1,1)+x((-j,y_j)+(1,1))\,|\,x\in\Q,x<0\}\subset\phi(P).
$$
In particular, taking $x=-t\in\Z$, it implies
$(t(j-1)-1,-(t(y_j+1)+1)\in\phi(P)$, i.e.,
$$
x_{t(y_j+1)+1}\ge t(j-1)-1\mbox{ for }t,j\in\Z_+\bs\{0\}.
\eqno(2.23)$$ Dividing (2.23) by $t(y_j+1)+1$ and taking
$\lim_{j\rar\infty}$, it gives $\b\ge\a^{-1}$. Thus $\b=\a^{-1}$.
\vs{-3pt}\par Assume that $\a=\OVER{q}{p}$ is a rational
number,\vs{-3pt} where $p,q\in\Z_+\bs\{0\}$ are coprime. By Lemma
2.9, there exists some $m\in\Z$ such that $(mp,-mq-1)=
(0,-1)+m(p,-q)\notin\phi(P)$. Say, $m>0$. Then again by Lemma 2.9,
$(0,0)+i(-mp,mq+1)\in\phi(P)$ for all $i\in\ool{0,\infty}$ since
$(0,0)\in\phi(P)$. \vs{-3pt} But taking $i>>0$, we have
$y_{imp}<i(mq+1)$ because
$\lim_{i\rar\infty}\OVER{y_{imp}}{imp}=\a=\OVER{q}{p}$.\vs{-3pt}
This contradicts the definition of $y_{imp}$ in (2.16).
\par
Thus $\a$ is not a rational number.
We define a well order $\DAYU$ on $\Z^2$ as follows:
$$
(i,j)\DAYU(k,l)\Rla(i-k,j-l)\DAYU(0,0)\mbox{ \ \ and \ \ }
(i,j)\DAYU(0,0)\Rla j>-i\a,
$$
i.e., $(i,j)\DAYU(0,0)$ if it
is located above the line $\{x(-1,\a)\,|\,x\in\R\}$
on the $Oxy$-plane.
\par
First we claim
\par
{\bf Claim 1.} For any $\e\!\in\!\R_+\bs\{0\}$, there exist $p,q\!\in\!\Z_+$
or $p,q\!\in\!\Z_-$
such that $0\!<\!q\!-\!p\a\!<\!\e$.
\par
It suffices to prove by induction on $n$ that there exist $p_n,q_n\in\Z$ such
that (such $p_n,q_n$ must have the same\vs{-9pt} sign)
$$\matrix{
|p_n\a-q_n|<(
{{\dis1}_{_{}}\over{}^{^{^{^{^{\!}}}}}{\dis2}{}^{^{^{^{^{\!}}}}}}
)^n\mbox{ for }n=1,2,...
\cr}
\vs{-9pt}\eqno(2.24)$$
Clearly, we can take $p_1\!=\!1$ and choose
$q_1\!\in\!\Z$ to satisfy $|\a\!-\!q_1|\!<\!
{{\dis1}_{_{}}\over{}^{^{^{^{^{\!}}}}}{\dis2}{}^{^{^{^{^{\!}}}}}}
$.
\vs{-6pt}Suppose (2.24) holds for $n$. Let $\a_n=|p_n\a-q_n|$. Choose
$r_n\in\Z$ such
that $|r_n\a_n-1|<
{{\dis1}_{_{}}\over{}^{^{^{^{^{\!}}}}}{\dis2}{}^{^{^{^{^{\!}}}}}}
\a_n$. \vs{-7pt}Let $p_{n+1}=r_np_n$,
$q_{n+1}=r_nq_n+1$, then we have
$|p_{n+1}\a-q_{n+1}|<
{{\dis1}_{_{}}\over{}^{^{^{^{^{\!}}}}}{\dis2}{}^{^{^{^{^{\!}}}}}}
\a_n<(
{{\dis1}_{_{}}\over{}^{^{^{^{^{\!}}}}}{\dis2}{}^{^{^{^{^{\!}}}}}}
)^{n+1}.\vs{-3pt}$ This proves the claim.\par Let
$\G_+=\{(i,j)\in\Z^2\,|\,(i,j)\DAYU(0,0)\}$. Then $\G_+$ is a
standard Borel subset of $\Z^2$ in the sense [10] that if
$(i,j),(k,l)\in\G_+$ such that $(i,j)\DAYU(k,l)$, then there
exists $n\in\Z_+$ such that $n(k,l)\DAYU(i,j)$.
\par
First assume that $\G_+\cap\phi(P)=\emptyset$, that is, there is
no weight located above the line \mbox{$\{x(-1,\a)\,|\,x\in\R\}$.}
Then it can be verified that $V$ is a Verma module defined in
[10], thus is not quasi-finite. To be self-contained, we simply
prove its non-quasi-finiteness as follows: For any $n>0$, using
Claim 1, we can choose $(p_i,q_i)\in\Z^2$ such that
$(0,-1)\XIYU(p_i,q_i)\XIYU(0,0)$ for $i=1,2,...,n$. Denote
$\mu_i=-(p_i,q_i)$, $\nu_i=-(0,-1)-\mu_i\DAYU(0,0)$. Recall
notation (2.2) that $L_\mu$ means $L_{ib_1+jb_2}$ if $\mu=(i,j)$,
we claim that $v_i=L_{-\nu_i}L_{-\mu_i}v_{\ssc\L}$, $i=1,2,...,n$,
all having the weight $\phi^{-1}(0,-1)$, are linear independent:
\vs{-7pt}suppose
$$
\SUM{\!i=1}{\,n}c_iL_{-\nu_i}L_{-\mu_i}v_{\ssc\L}=0\mbox{ for some
}c_i\in\C. \vs{-4pt}\eqno(2.25)$$ Say, with respect to the
ordering $\DAYU$, $\nu_1$ is the largest elements among all
$\mu_i,\nu_i$ with $c_i\ne0$. Claim 1 shows that one can choose
$\g\in\Z^2$ such that all $\mu_i,\nu_i\XIYU\g$ except that
$\g\XIYU\nu_1$. Applying $L_\g$ to (2.25), using that
$\G_+\cap\phi(P)=\emptyset$, we obtain an equation which has the
\vs{-5pt}form
$$
c_1L_{\g-\nu_1}L_{-\mu_1}v_{\ssc\L}+
c'L_{\g-\nu_1-\mu_1}v_{\ssc\L}=0\mbox{ for some }
c'\in\C.
\vs{-4pt}\eqno(2.26)$$
By applying $L_{\nu_1+\mu_1-\g}$ and $L_{\mu_1}L_{\nu_1-\g}$ to (2.26), we
obtain $v_{\ssc\L}=0$, a contradiction.
\par
Next assume that there exists $(i,j)\in\G_+\cap\phi(P)$. Then for
any $(k,l)\XIYU(i,j)$, we must have $(k,l)\in\phi(P)$, otherwise
by Lemma 2.9 we would have $(t(i-k)+i,t(j-l)+j)=(i,j)+t(i-k,j-l)
\in\phi(P)$ for all $t\ge0$, which gives the following
contradiction: when $t>>0$, either $(t(i-k)+i,t(j-l)+j)\ge0$, or
$t(i-k)+i<0$ but $y_{t(k-i)-i}<t(j-l)+j$ by definition of $\a$ in
(2.21), or $t(j-l)+j<0$ but $x_{t(l-j)-j}<t(i-k)+i$ by definition
of $\b$.
\par
We have $i<0<j$ \vs{-3pt}or $i>0>j$ by (2.6). Say $i<0$, and
rewrite $(i,j)$ as $(-i,j)$. For a given $n>0$, let
$\e=\OVER{1}{n}(j-i\a)>0$. \vs{-3pt}Choose $p,q$ to be as in Claim
1. Then we obtain $(0,0)\XIYU(-p,q)$ and $(-np,nq)\XIYU(-i,j)$.
Note that since $q>\a p$, when $m>>0$ we have $y_{mp}<mq$ if
$p,q>0$ or $x_{-mq}<-mp$ if $p,q<0$, i.e.,
$(-mp,mq)\notin\phi(P)$. So we can let $m\ge n$ be the largest
integer such that $(-mp,mq)\in\phi(P)$. Let $v_\l$ be a nonzero
weight vector with weight $\l=\phi^{-1}(-mp,mq)$, then $v_\l$
generates a nontrivial highest weight $\Vir[-pb_1+qb_2]$-submodule
$V'$ of $V$, where $\Vir[-pb_1+qb_2]={\rm
span}\{L'_i=L_{-ip,iq}\,|\,i\in\Z\}$ (if $V'$ is trivial, then it
is a trivial module over a rank two Virasoro algebra
$\Vir[-pb_1+qb_2,kb_1+kb_2]$ which is generated by
$\{L_{\pm(-pb_1+qb_2)},L_{kb_1+kb_2}\}$ for some $k>>0$,
contradicting Lemma 2.6), such that
$\L=\phi^{-1}(0,0)=\l-m(-pb_1+qb_2)$ is a weight with weight
multiplicities at least $m$ by the well-known property of a
nontrivial highest weight module over the Virasoro algebra. Since
$m\ge n$ and $n$ is arbitrary, we obtain that the weight
multiplicity of $\L$ is infinite and so $V$ is not quasi-finite.
This proves that a nontrivial irreducible GHW $\HVir$-module does
not exist and thus we obtain Theorem 1.1 for the special case
$n=2$.
\par\
\vs{-5pt}\par
\ni{\bf
             3. \ Proofs of the Theorems in General Case
}
\vs{-1pt}\par\ni
Now we can complete the proofs of Theorems 1.1, 1.2 as follows.
\par\ni
{\it Proof of Theorem 1.1.} Let $V$ be a not-finitely-dense
Harish-Chandra module over the higher rank Virasoro algebra
$\HVir$, where $M$ has rank $n\ge2$. Suppose $V$ is not a module
of the intermediate series. Then by Lemma 2.1(2), there exists a
GHW vector $v_{\ssc\L}$ and some $\Z$-basis $B=\{b_1,...,b_n\}$ of
$M$ such that (2.1) holds. Assume that $V$ is nontrivial, then
there exists at least some $i,\,1\le i\le n$ such that
$L_{-b_i}v_{\ssc\L}\ne0$. Choose any $j\ne i$, then $v_{\ssc\L}$
generates a nontrivial GHW $\Vir[b_i,b_j]$-submodule $V'$ of $V$
such that $v_{\ssc\L}$ is a nontrivial GHW vector. Therefore at
least one of composition factors of the $\Vir[b_i,b_j]$-module
$V'$ is a nontrivial irreducible GHW $\Vir[b_i,b_j]$-module,
contradicting the result we obtained for the case $n=2$.
\qed\par\ni {\it Proof of Theorem 1.2.} Suppose $V$ is a
not-finitely-dense indecomposable quasi-finite weight modules over
$W(m,\G)$ or over $\wh W(1,\G)$, where $m\ge2$ or there exists a
group injection $\Z\times\Z\rar\G$, i.e., there exists $\Z$-linear
independent elements $\a_1,\a_2$ of $\G$. Below we shall only
consider the case $W(m,\G)$ as $\wh W(1,\G)$ is the central
extension version of $W(1,\G)$. Using notation as in Introduction,
$W(m,\G)={\rm span}\{ x^\a\ptl_i\,|\,\a\in\G,i=1,2,...,m\}$ such
\vs{-4pt}that
$$
[x^\a\ptl_i,x^\b\ptl_j]=x^{\a+\b}(\pi_i(\b)\ptl_j-\pi_j(\a)\ptl_i)
\mbox{ for }\a,\b\in\G,\,i,j=1,2,...,m,
\vs{-4pt}$$
where $\pi_i:\C^m\rar C$ is the natural projection and $\G\subset\C^m$ is
a nondegenerate subgroup.
Let $\DD={\rm span}\{\ptl_i\,|\,i=1,2,...,m\}$ be the Cartan subalgebra
of $W(m,\G)$. Then $V$ has the weight space
\vs{-4pt}decomposition:
$$
V=\OPLUS{\!\!\!\!\l\in\DD^*}V_\l,\ \ V_\l=\{v\in
V\,|\,\ptl{\ssc\,} v=\l(\ptl)v \mbox{ for }\ptl\in\DD\},
\vs{-4pt}$$ such that all $V_\l$ are finite dimensional,
where $\DD^*$ is the dual space of $\DD$. We regard $\G$ as an
additive subgroup of $\DD^*$ by defining $\a(\ptl_i)=\pi_i(\a)$
for $\a\in\G,i=1,2,...,m$. For any $\mu\in\DD^*$, we let
$V[\mu]=\oplus_{\l\in\mu+\G}V_\l$. Then $V[\mu]$ is a submodule of
$V$ such that $V$ is a direct sum of different $V[\mu]$. Thus
$V=V[\mu]$ for some $\mu\in\DD^*$. If $0$ is a weight of $V$, we
take $\mu=0$.
\par
Suppose $\l$ is a weight. Define
a linear transformation $T_\l:V_\l\rar
V'=V_\mu+V_{\mu+\a_1}+V_{\mu+\a_2}$ as follows: If $\l=\mu$,
we take $T_\mu=1_{V_\mu}$ to be the
embedding. Suppose $\l\ne\mu$. Then also $\l\ne0$.
Let $\a=\mu-\l\in\G$. We choose $\b=\a+\a_1$ or
$\a+\a_2$ such that $\a,\b$ are $\Z$-linear independent. From linear algebra,
there exists $\ptl\in\DD$ with $b_1=\a(\ptl),b_2=\b(\ptl)\in\C$ being $\Z$-linear
independent and $\l(\ptl)\in\C\bs\{0\}$.
\vs{-4pt}Define
$$
T_\l:V_\l\rar V'\mbox{ such that }
T_\l(v)=(x^\a\ptl){\ssc\,}v+(x^\b\ptl){\ssc\,}v\in V' \mbox{ for
}v\in V_\l. \vs{-4pt}\eqno(3.1)$$ We claim that $T_\l$ is an
injection. Suppose $v\ne0$ such that $T_\l(v)=0$. Let
$\G'=\{m\a+n\b\,|\,m,n\in\Z\}$ be a rank two subgroup of $\G$.
Then $\ol{\Vir}[\G']={\rm span}\{L_\g=x^\g\ptl\,|\,\g\in\G'\}$ is
a rank two (centerless) Virasoro algebra isomorphic to the
centerless version of $\Vir[b_1,b_2]$. Take a $\Z$-basis
$\{\g_1=\a+\b, \g_2=\a+2\b\}$ of $\G'$. Then from $L_\a v=L_\b
v=0$, we obtain $L_{m\g_1+n\g_2}v=0$ for all $n,m\ge0$ such that
$(m,n)\ne(0,0)$. Thus $v$ is a GHW vector. Since
$L_0v=\ptl{\ssc\,}v=\l(\ptl)v\ne0$, $v$ generates a nontrivial GHW
$\ol{\Vir}[\G']$-submodule of $V$, a contradiction with the result
we obtained for the case $n=2$. Thus $T_\l$ is an injection.
Therefore we have ${\rm dim\ssc\,}V_\l\le {\rm dim\ssc\,}V'$ for
all $\l\in\DD^*$, i.e., $V$ is a uniformly bounded module. \par
Now assume that ${\rm dim\,}V_{\mu+\a}\ne{\rm dim\,}V_{\mu+\b}$
for some $\a,\b\in\G$ such that $\a\ne\b$ and
$\mu+\a\ne0\ne\mu+\b$. Choose $\ptl\in\DD$ such that
$(\mu+\a)(\ptl),(\b-\a)(\ptl),(\mu+\b)(\ptl)\ne0$. Form a rank one
centerless Virasoro subalgebra $\ol{\Vir}[\b-\a]={\rm
span}\{x^{i(\b-\a)}\ptl\,|\,i\in\Z\}$ and let
$V'=\oplus_{i\in\Z}V_{\mu+\a+i(\b-\a)}$ be a
$\ol{\Vir}[\b-\a]$-submodule of $V$. Then regarding $V'$ as a
$\ol{\Vir}[\b-\a]$-module, each subspace $V_{\mu+\a+i(\b-\a)}$ is
exactly the weight space of weight $(\mu+\a+i(\b-\a))(\ptl)$ with
respect to $\ptl\in\ol{\Vir}[\b-\a]$ (since for different $i$,
$(\mu+\a+i(\b-\a))(\ptl)$ have different values, and note that
$\C\ptl$ is the Cartan subalgebra of $\ol{\Vir}[\b-\a]$). Thus
$V'$ is also a uniformly bounded module over $\ol{\Vir}[\b-\a]$.
Thus $V'$ must only have a finite composition factors (and each
nontrivial composition factor is a module of the intermediate
series), and we see that all the nonzero weights (with respect to
$\ptl$) of $V'$ have the same multiplicity. This is a
contradiction with ${\rm dim\,}V_{\mu+\a}\ne{\rm
dim\,}V_{\mu+\b}$.\qed\par\ \vs{-5pt}\par%
\ni{\bf
                        References
} \small\vs{-1pt} \par\ni\hi3.5ex\ha1\lineskip 4.5pt %
[1] B.~N.~Allison, S.~Azam, S.~Berman, Y.~Gao and A.~Pianzola,
Extended affine Lie algebras and their root systems, {\it
Mem.~Amer.~Math.~Soc.} {\bf 126}, No. 603, 1997.
\par\ni\hi3.5ex\ha1 
[2] V.~Chari and A.~Pressley,
 Unitary representations of the Virasoro algebra and a conjecture of Kac,
 {\it Compositio Math.} {\bf67} (1988), 315-342.
 \par\ni\hi3.5ex\ha1
[3] D.~Friedan, Z.~Qiu and S.~Shenker, Conformal invariance,
unitarity and two-dimensional critical exponents, in {\it Vertex
Operators in Mathematics and Physics} (J.~Lepowsky, S. Mandelstam
and I.~M.~Singer, Eds.), Springer, New York-Berlin, 1985.
 \par\ni\hi3.5ex\ha1 
[4] I.~M.~Gelfand and D.~B.~Fuchs, Cohomologies of the Lie algebra
of vector
 fields on the circle,
 {\it Funct.~Anal.~Appl.} {\bf2} (1968), 92-39 (English translation 114-126).
 \par\ni\hi3.5ex\ha1
[5] P.~Goddard and D.~Oliver, Kac-Moody and Virasoro algebras in
relation to quantum physics, {\it Internat.~J.~Mod.~Phys.} {\bf
A1} (1986), 303-414.
 \par\ni\hi3.5ex\ha1 
[6] V.~G.~Kac, Some problems of infinite-dimensional Lie algebras
and
 their representations, {\it Lecture Notes in Mathematics} {\bf933}, 117-126,
 Springer, 1982.
 \par\ni\hi3.5ex\ha1 
[7] V.~G.~Kac, {\it Infinite Dimensional Lie Algebras}, 3rd ed.,
  Combridge Univ. Press, 1990.
  \par\ni\hi3.5ex\ha1 
[8] C.~Martin and A.~Piard, Indecomposable modules over the
Virasoro Lie
 algebra and a conjecture of V Kac,
{\it Comm. Math. Phys.} {\bf137} (1991), 109-132.
 \par\ni\hi3.5ex\ha1 
[9] O.~Mathieu, Classification of Harish-Chandra modules over the
 Virasoro Lie algebra, {\it Invent. Math.} {\bf107} (1992), 225-234.
 \par\ni\hi4ex\ha1 
[10] V.~Mazorchuk, Verma modules over generalized Witt algebras,
{\it Compositio Math.} {\bf115} (1999), 21-35.
 \par\ni\hi4ex\ha1 
[11] V.~Mazorchuk, Classification of simple Harish-Chandra modules
over
  $\Q$-Virasoro algebra, {\it Math. Nachr.} {\bf209} (2000), 171-177.
  \par\ni\hi4ex\ha1 
[12] A.~Neveu and J.~H.~Schwarz, Factorizable dual model of pions,
 {\it Nucl. Phys. B} {\bf31} (1971), 86-112.
 \par\ni\hi4ex\ha1 
[13] J.~Patera and H.~Zassenhaus, The higher rank Virasoro
algebras, {\it Comm. Math. Phys.} {\bf136} (1991), 1-14.
 \par\ni\hi4ex\ha1 
[14] P.~Ramond, Dual theory of free fermions, {\it Phys. Rev. D}
 {\bf3} (1971), 2451-2418.
 \par\ni\hi4ex\ha1
[15] C.~Rebbi, Dual models and relativistic quantum strings, {\it
Phys.~Rep.} {\bf12} (1974), 1-73.
 \par\ni\hi4ex\ha1 
[16] Y.~Su, A classification of indecomposable $sl_2(\C)$-modules
 and a conjecture of Kac on irreducible modules over the Virasoro
 algebra, {\it J. Alg.} {\bf161} (1993), 33-46.
 \par\ni\hi4ex\ha1 
[17] Y.~Su, Harish-Chandra modules of the intermediate series over
the high
 rank Virasoro algebras and high rank super-Virasoro algebras,
 {\it J. Math. Phys.} {\bf35} (1994), 2013-2023.
 \par\ni\hi4ex\ha1 
[18] Y.~Su, Classification of Harish-Chandra modules over the
 super-Virasoro algebras, {\it Comm. Alg.} {\bf23} (1995), 3653-3675.
 \par\ni\hi4ex\ha1 
[19] Y.~Su, Simple modules over the high rank Virasoro algebras,
{\it Comm.~Alg.} {\bf29} (2001), 2067-2080.
 \par\ni\hi4ex\ha1 
[20] Y.~Su, Simple modules over the higher rank super-Virasoro
algebras,
 {\it Lett. Math. Phys.} {\bf53} (2000), 263-272.
 \par\ni\hi4ex\ha1 
[21] Y.~Su, On indecomposable modules over the Virasoro algebra,
 {\it Science in China A} {\bf44} (2001), 980-983 (also available at Math.QA/0012013).
 \par\ni\hi4ex\ha1 
[22] Y.~Su and K.~Zhao, Generalized Virasoro and super-Virasoro
algebras and modules of the intermediate series, {\it J.~Alg.}
{\bf252} (2002), 1-19.
 \par\ni\hi4ex\ha1 
[23] X.~Xu, New generalized simple Lie algebras of Cartan type
over a field
  with characteristic 0, {\it J. Alg.} {\bf224} (2000), 23-58.
\end{document}